\documentclass[conference,9pt]{IEEEtran}
\ifCLASSINFOpdf
\else
\fi
\IEEEoverridecommandlockouts

\usepackage[cmex10]{amsmath}
\usepackage{amsfonts,amssymb}
\usepackage{amsthm}
\usepackage{tikz}
\usepackage{pgfplots}
\usepackage{array}

\ifCLASSOPTIONcompsoc
 \usepackage[caption=false,font=normalsize,labelfont=sf,textfont=sf]{subfig}
\else
 \usepackage[caption=false,font=footnotesize]{subfig}
\fi
\theoremstyle{plain}
\newtheorem{thm}{Theorem}
\newtheorem{lem}[thm]{Lemma}
\newtheorem{col}[thm]{Corollary}

\theoremstyle{definition}
\newtheorem{assum}{Assumption}
\newtheorem{dfn}{Definition}

\theoremstyle{remark}
\newtheorem{rem}{Remark}

\usepackage{mydefs}
\usepackage{color,soul}
\usepackage{tabu,multirow}
\usepackage{cite}

\begin{document}
\title{A Unified Contraction Analysis of a Class of Distributed Algorithms for Composite Optimization\vspace{-1.1cm}}
 \author{Jinming Xu, Ying Sun, Ye Tian, and Gesualdo Scutari\thanks{School of Industrial Engineering, Purdue University, West-Lafayette, IN, USA. Emails: \texttt{<xu1269, sun578, tian110, gscutari>} \texttt{@purdue.edu.} This work has been supported by the  USA NSF
Grants CIF 1632599 and CIF 1564044; and the ARO  Grant
W911NF1810238.}}

\maketitle

\begin{abstract}
We study distributed composite  optimization over networks: agents minimize  the sum of a smooth (strongly) convex  function--the agents' sum-utility--plus a nonsmooth (extended-valued) convex one. We propose a  general  algorithmic framework for such a class of problems and provide a   unified convergence analysis leveraging   the theory of operator splitting. Our results unify several approaches proposed in the literature of distributed optimization for special instances of our formulation.  Distinguishing features of our scheme are: (i)  when the agents' functions are strongly convex, the algorithm converges at a {\it linear} rate,  whose dependencies on the agents' functions   and the network topology are {\it decoupled}, matching the typical rates of   centralized optimization; 
(ii)  the step-size  does not depend on the network parameters but only on the optimization ones; and (iii) the algorithm can adjust the ratio between the number of communications and computations to achieve the {\it same} rate of the centralized proximal gradient scheme (in terms of computations). This is the first time that a distributed algorithm applicable to {\it composite} optimization enjoys such properties.
\vspace{-0.3cm}
\end{abstract}


%
\IEEEpeerreviewmaketitle

 \section{Introduction}\vspace{-0.2cm}
We study distributed multi-agent  optimization over networks, modeled as undirected static graphs. Agents aim at solving \vspace{-0.2cm}
\begin{equation}\label{prob:dop_nonsmooth_same}
\min_{x\in\mathbb{R}^{d}} F(x)+G(x),\quad F(x)\triangleq {\frac{1}{m}}\sum_{i=1}^m f_i(x),\tag{P}\vspace{-0.1cm}
\end{equation}
where $f_i:\mathbb{R}^d\to \mathbb{R}$ is the cost-function of agent $i$, assumed to be smooth, (strongly) convex and known only to the agent; and $G:\mathbb{R}^d\to \mathbb{R}\cup\{-\infty, \infty\}$ is a  nonsmooth, convex (extended-value) function, which can be used to enforce shared constraints  or some specific structure on the solution, such as sparsity.\\
\indent Our focus is on the design of   distributed algorithms for Problem \eqref{prob:dop_nonsmooth_same} that provably converge at a  {\it linear}   rate.
When $G=0$,  several distributed schemes have been proposed in the literature enjoying such a property; examples include EXTRA~\cite{shi2015extra}, AugDGM~\cite{xu2015augmented}, NEXT~\cite{di2016next}, SONATA \cite{YingMAPR,sun2019convergence}, DIGing~\cite{nedich2016achieving},   NIDS~\cite{li2017decentralized},   Exact Diffusion~\cite{yuan2018exact_p1},  MSDA~\cite{scaman17optimal}, and the distributed algorithms in \cite{qu2017harnessing},\cite{jakovetic2018unification}, and \cite{mansoori2019general}.  When $G\neq 0$ results are scarce; to our knowledge, the only two schemes available in the literature  achieving linear rate for \eqref{prob:dop_nonsmooth_same} are SONATA \cite{sun2019convergence} and the distributed proximal gradient algorithm  \cite{alghunaim2019linearly}. 
The aforementioned algorithms apparently look different; no unified convergence analysis can be inferred; and, in most of the cases,  step-size bounds and convergence rate seem quite conservative. This naturally suggests the following two questions: 
\vspace{-0.2cm}\begin{itemize}
	\item[{\bf (Q1)}] Can one unify the design and analysis of distributed algorithms in the setting \eqref{prob:dop_nonsmooth_same}?
	\vspace{-0.1cm}\item[{\bf (Q2)}]  Can one match the linear convergence rate of the centralized proximal-gradient algorithm applied to   \eqref{prob:dop_nonsmooth_same}?
\vspace{-0.2cm}
\end{itemize}
Recent efforts toward a better understanding of the taxonomy of  distributed algorithms (question Q1) are the following: \cite{jakovetic2018unification} provides a connection between EXTRA and DIGing; \cite{Scoy-canonical18} provides a canonical representation of some of the distributed algorithms above--NIDS and Exact-Diffusion are proved to be equivalent; and \cite{Sundararajan17} provide an automatic (numerical) procedure to prove linear rate of some classes of distributed algorithms.  These efforts model only first order algorithms applicable to Problem \eqref{prob:dop_nonsmooth_same} {\it with $G=0$} and employing a {\it single} round of communication and gradient computation. Because of that, in general, they cannot achieve the rate of the centralized gradient algorithm (addressing thus Q2). Works partially addressing Q2 are the following:   MSDA~\cite{scaman17optimal} uses multiple communication steps to achieve the lower complexity bound of  \eqref{prob:dop_nonsmooth_same} when $G=0$; and
the algorithms in \cite{van2019distributed} and \cite{li2017decentralized} achieve linear rate and can adjust the number of communications performed at each iteration to  match the rate of the centralized gradient descent. 
However it is not clear how to extend  (if possible) these methods and their convergence analysis to the more general composite (i.e., $G\neq 0$) setting \eqref{prob:dop_nonsmooth_same}. \\\indent This paper aims at addressing Q1 and Q2 
in the general setting \eqref{prob:dop_nonsmooth_same}. Our major contributions are the following: 1) We propose a general primal-dual  distributed  algorithmic framework   that subsumes several existing ATC- and CTA-based distributed algorithms; 2) A sharp linear convergence rate is proved (when $G\neq 0$) developing  an operator contraction-based analysis. 
By product, our convergence results apply also to the algorithms in \cite{shi2015extra,li2017decentralized,yuan2018exact_p1,qu2017harnessing,di2016next,xu2015augmented}, which so far have been studied in isolation; 3) For ATC forms    of our schemes, the dependencies of the linear rate  on the agents' functions   and the network topology are {\it decoupled}, matching the typical rates for the centralized optimization and the consensus averaging. This is a major departure from  existing analyses,  which do not show such a clear separation, and complements the results in \cite{li2017decentralized} applicable only to smooth instances of  \eqref{prob:dop_nonsmooth_same}. Furthermore,  convergence is established under a proper choice of the step-size, whose upper bound does {\emph not} depend on the network parameters but only on the optimization ones (Lipschitz constants of the gradients and strongly convexity constants); and 4) The proposed scheme can naturally adjusts the ratio between the number of communications and computations to achieve the {\it same} rate of the centralized proximal gradient scheme (in terms of computations). Chebyshev acceleration can also be  employed to significantly reduce the number of communication steps per computation. Because of space limitation, all the proofs are available as supporting material in the technical report \cite{XuSunScutariJ}.\\
\noindent {\bf Notations}:   $\mathbb{N}_+$ is the  set of positive integer numbers; $\mathbb{S}^m$ is the set of $\mathbb{R}^{m\times m}$ symmetric matrices   while  $\mathbb{S}_+^m$ (resp.  $\mathbb{S}_{++}^m$) is the set of   positive semidefinite  (resp. definite) matrices in $\mathbb{S}^m$.  $\mathbb{P}_K$  denotes  the set of (real) monic polynomials of order $K$. {Unless otherwise indicated,  column vectors are denoted by lower-case letters while upper-case letters are used for  matrices (with the exception of $L$ in Assumption 1 to conform with conventional notation). 
The symbols $1_m$ and $0_m$ denote the $m$-length column vectors of all ones and all zeros, respectively. The $\mathbf{0}_m$ denotes the $m\times m$ zero matrix;   $I_m$ denotes the  identity matrix in $\mathbb{R}^{m\times m}$; $\J\triangleq 1_m\,1_m^\top/m$ is the projection matrix onto $1_m$.}   With a slight abuse of notation,  $I$ will denote either the identity matrix or the identity operator on the space under consideration.   We use $\Null{\cdot}$ [resp. $\Span{\cdot}$] to denote the null space (resp. range space) of the matrix argument.   
For any  $X,Y\in\mathbb{R}^{m\times d}$, let $\innprod{X}{Y}\triangleq\text{trace}(X^\top Y)$ while we write  $\norm{X}$ for $\norm{X}_F$;   the same notation is used for vectors, treated as special cases.
Given ${G}\in  \mathbb{S}_+^n$,    $\innprod{X}{Y}_\G\triangleq\innprod{\G\x}{\y}$  {and} $\norm{\x}_\G\triangleq\sqrt{\innprod{\x}{\x}_{\G}}.$ The eigenvalues of a symmetric matrix ${A}\in \mathbb{R}^{m\times m}$ are denoted by  $\lambda_i(\A)$, $i=1,\ldots, m$, and arranged in increasing order. For $x\in \mathbb{R},$ we denote $x_+ = \max(x,0).$
\vspace{-0.3cm}


\section{Problem Statement}
\vspace{-0.2cm}
We study Problem~\eqref{prob:dop_nonsmooth_same} under the following  assumption.
\begin{assum}\label{assum:Lipschiz_gradient}
Each local cost function $f_i:\mathbb{R}^d\rightarrow\mathbb{R}$ is  $\mu$-strongly convex  and $L$-smooth; and $G:\mathbb{R}^d\rightarrow\mathbb{R}\cup\{\pm\infty\}$ is proper, closed and convex. Define $\kappa\triangleq L/\mu$.\vspace{-0.4cm}
\end{assum} 

Note that Assumption \ref{assum:Lipschiz_gradient} also accounts for the case where $f_i$ is convex and $G$ is $\mu$-strongly convex.\\
\noindent {\bf Network model:} Agents are embedded in a network, modeled as  an undirected, static graph $\Gh=(\Vx,\Eg)$, where $\Vx$ is the set of nodes (agents) and $\{i,j\}\in\Eg$ if there is an edge (communication link) between node $i$ and $j$. 
We make the blanket assumption that $\Gh$ is connected.  We introduce the following matrices associated with $\mathcal{G}$, which will be used to build the proposed distributed algorithms. \vspace{-0.1cm}
\begin{dfn}[Gossip matrix]
\label{dfn:weight_matrix}
A matrix $\W\triangleq [W_{ij}]\in \mathbb{R}^{m\times m}$ is said to be compliant to the graph $\Gh=(\Vx,\Eg)$ if  $W_{ij} \neq 0$ for  $\{i,j\}\in\Eg$, and $W_{ij}=0$ otherwise. The set of such matrices is denoted by $\mathcal{W}_\mathcal{G}$.
\end{dfn}
\begin{dfn}[$K$-hop gossip matrix]
\label{dfn:k-hop_weight_matrix}
Given $K\in\mathbb{N}_+,$ a matrix $ {W}'\in\mathbb{R}^{m\times m}$ is said to be a $K$-hop gossip matrix associated to $\Gh=(\Vx,\Eg)$  if $ {W}'=P_K(\W)$, for some $\W\in\mathcal{W}_\Gh$, where   $P_K(\cdot)\in\mathbb{P}_K$.\vspace{-0.1cm}
\end{dfn}
\vspace{-0.2cm}
Note that, if  $\W\in \mathcal{W}_\Gh$,  using $W_{ij}$ to   linearly combine information between   agent $i$ and $j$ corresponds to performing a single communication  between  the two agents ($i$ and $j$ are immediate neighbors).   Using  a $K$-hop matrix ${W}'=P_K(\W)$  requires instead $K$ consecutive rounds of communications among immediate neighbors for the   aforementioned weighting process to be implemented in a distributed way (note that the zero-pattern of $ {W}'$ is in general not compliant with $\mathcal{G}$).
   $K$-hop weight matrices are crucial to employ acceleration of the communication step, which will be a key ingredient to exploit the tradeoff between  communications and computations (cf.~Sec.~\ref{sec:tradeoff}).\\
\noindent{\textbf{A saddle-point reformulation:}} Our path to design distributed solution methods for \eqref{prob:dop_nonsmooth_same} is to solve a   saddle-point reformulation of  \eqref{prob:dop_nonsmooth_same}   via general proximal splitting algorithms that are implementable over  $\mathcal{G}$.   Following a standard path in the literature,
we introduce local copies $x_i\in \mathbb{R}^d$ (the $i$-th one is owned by agent $i$) of  $x$  and functions \vspace{-0.2cm}\begin{equation}
	f(\x)\triangleq\sum_{i=1}^mf_i(x_i)\quad \text{and}\quad g(\x)\triangleq\sum_{i=1}^m G(x_i),\vspace{-0.1cm}
\end{equation} with $\x\triangleq [x_1,\ldots, x_m]^\top\in\mathbb{R}^{m\times d}$;
  \eqref{prob:dop_nonsmooth_same} can be then rewritten as\vspace{-0.1cm}
\begin{equation}\label{prob:dop_nonsmooth_same_augmented}
\min_{\x\in\mathbb{R}^{m\times d}} f(\x)+g(\x),~{\text{s.t.}}~  \sqrt{{C}}\x=\zeros,\vspace{-0.1cm}
\end{equation}
where $C$ satisfies the following assumption:
\begin{assum}\label{assum:cond_C}
$C \in \mathbb{S}^m_+$ and $\Null{{C}}=\Span{1}$.
\end{assum}
Under this condition,  the constraint $\sqrt{ {C}}\x=\zeros$   enforces a consensus among $x_i$'s and thus  \eqref{prob:dop_nonsmooth_same_augmented} is equivalent to   \eqref{prob:dop_nonsmooth_same}. 
\\ \indent 
In the setting above, \eqref{prob:dop_nonsmooth_same_augmented} is equivalent to its KKT conditions: there exists $\x^\star \in  \mathcal{S}_{\text{KKT}}$, where $\mathcal{S}_{\text{KKT}}$ is defined as\vspace{-0.1cm}
\begin{align}\label{eq:kkt_conditions}
 \mathcal{S}_{\texttt{KKT}}  \triangleq & \left\{   \x \in  \mathbb{R}^{m\times d} \, \big\vert\, \exists\, \y \in  \mathbb{R}^{m\times d} \text{ such that} \right.\nonumber\\
& \left.  \sqrt{{C}}\x=\zeros, \quad  \nabla f(X)+  \sqrt{{C}}\y \in-\partial g(\x)\right\},
\end{align}
where $\nabla f(\x)\triangleq [\nabla f_1(x_1),\nabla f_2(x_2),...,\nabla f_m(x_m)]^\top$
and  $\partial g(\x)$ denotes the subdifferential  of $g$ at $\x$. We have the following. 

\begin{lem}\label{lemma_eq_KKT} Consider   Problem \eqref{prob:dop_nonsmooth_same}  under Assumptions 1 and 2; $x^\star\in \mathbb{R}^d$ is an optimal solution of \eqref{prob:dop_nonsmooth_same} if and only if  ${1}_m x^{\star\top} \in \mathcal{S}_{\texttt{KKT}}$.  \vspace{-0.3cm}
\end{lem}

Building on Lemma \ref{lemma_eq_KKT}, 
in the next section, we propose a general distributed algorithm for \eqref{prob:dop_nonsmooth_same}  based on  a suitably defined operator splitting solving  the KKT system   \eqref{eq:kkt_conditions}. 
 \vspace{-0.2cm}

\section{A General Primal-Dual Proximal Algorithm}\label{sec:alg-des}\vspace{-0.1cm}

The proposed general primal-dual proximal algorithm reads
\begin{subequations}\label{alg:g-pd-ATC}
\begin{align}
\x^k&=\prox{\gamma g}{\z^k}, \label{alg:g-pd-ATC_z}\\
\z^{k+1}&=\A\x^k-\gamma\B\nabla f(\x^k)-\y^k, \label{alg:g-pd-ATC_x}\\
\y^{k+1}&=\y^k+{\C}\z^{k+1}, \label{alg:g-pd-ATC_y}
\end{align}
\end{subequations}
with  $\z^0\in \mathbb{R}^{m\times d}$ and $\y^0\in {\texttt{span}}({C})$. In \eqref{alg:g-pd-ATC_z},  $\prox{\gamma g}{\x}\triangleq\text{arg}\min_{\y}g(\y)+\frac{1}{2\gamma}\norm{\x-\y}^2$ is the standard proximal operator, which accounts for the nonsmooth term. Eq. \eqref{alg:g-pd-ATC_z} represents the update of the primal variables, where  $\A,\B\in \mathbb{R}^{m\times m}$  are suitably chosen weight matrices, and $\gamma>0$ is the step-size.   Finally,  \eqref{alg:g-pd-ATC_y} represents the update of  the dual variables.
Note that there is no loss of generality in initializing $\y^0\in {\texttt{span}}({C})$, as any  $\y$ in (\ref{eq:kkt_conditions}) is so  (unless all the $f_i$ share a common minimizer). 
\\\indent
Define the  set
$\mathcal{S}_{\texttt{Fix}}  \triangleq   \big\{   \x \in  \mathbb{R}^{m\times d} \, \big\vert\, \C\x=0 \text{ and }
  1^\top(I-\A)\x+\gamma \, 1^\top\B\nabla f(\x)\in - \gamma\,  1^\top\partial g(\x) \big\}$. It is not difficult to check that any fixed point $(X^\star, Z^\star, Y^\star)$ of  Algorithm \eqref{alg:g-pd-ATC} satisfies $X^\star \in \mathcal{S}_{\texttt{Fix}}$.
The following are {\it necessary} and sufficient conditions on  $\A$ and  $\B$   for  $X^\star \in \mathcal{S}_{\texttt{Fix}}$   to be the solution of \eqref{prob:dop_nonsmooth_same_augmented}.  \vspace{-0.1cm}
\begin{assum}\label{assum:cond_A_B}
The weight matrices $\A,\B \in \mathbb{R}^{m\times m}$ satisfy:
$1^\top\A \,1=m$,  and  $1^\top\B = 1^\top$.
\end{assum}
\begin{lem}[\!\cite{XuSunScutariJ}]
Under Assumption~\ref{assum:cond_C},  $\mathcal{S}_\texttt{KKT}=\mathcal{S}_\texttt{Fix}$ if and only if $\A,\B$ satisfy Assumption~\ref{assum:cond_A_B}.
\vspace{-0.3cm}
\end{lem}

\subsection{Connections with existing distributed algorithms} \label{Sec:special-cases}\vspace{-0.1cm}
Algorithm~\eqref{alg:g-pd-ATC} contains a gamut of distributed (and centralized) schemes, corresponding to different choices of the weight matrices $\A,\B,$ and $\C$; any   $A,B,C\in \mathcal{W}_{\mathcal{G}}$ leads to  distributed  implementations. The use of general matrices $\A$ and $\B$ (rather the more classical choices $\A=\B$ or $\B=I$) permits to model for the first time in a unified   algorithmic framework both ATC- and CTA-based  updates; this includes several existing distributed algorithms proposed for special cases of (P), as  briefly  discussed next; see \cite{XuSunScutariJ} for more examples. 
  Rewrite Algorithm~\eqref{alg:g-pd-ATC} in the following equivalent form: 
\begin{equation}\label{eq:g-pd-ATC_eliminate_y}
\z^{k+2}=(\I-\C)\z^{k+1}+\A(\x^{k+1}-\x^k)-\gamma\B(\nabla f(\x^{k+1})-\nabla f(\x^k)).\vspace{-0.2cm}
\end{equation}
When $G=0$, the above update reduces to
\begin{equation}\label{eq:g-pd-ATC_eliminate_y_G=0}
\x^{k+2}=(\I-\C+\A)\x^{k+1}-\A\x^k-\gamma\B(\nabla f(\x^{k+1})-\nabla f(\x^k)).\vspace{-0.1cm}
\end{equation}
It is not difficult to check that   the schemes in \cite{shi2015extra,li2017decentralized,yuan2018exact_p1,di2016next,xu2015augmented,nedich2016achieving,qu2017harnessing,jakovetic2018unification,mansoori2019general,alghunaim2019linearly} are all special cases of \eqref{eq:g-pd-ATC_eliminate_y} or \eqref{eq:g-pd-ATC_eliminate_y_G=0} and thus of Algorithm \eqref{alg:g-pd-ATC}--Table \ref{tab:connections_to_existing_algs} shows the proper parameter setting to establish the equivalence, where $\W\in  \mathcal{W}_\mathcal{G}$ is the weight matrix used in the target distributed algorithms, see \cite{XuSunScutariJ}  for more details.\vspace{-0.3cm}
\begin{table}
\renewcommand{\arraystretch}{1.2}
\scriptsize\centering
\setlength\tabcolsep{2pt}
\begin{tabu}{lc}
{\small \bf Algorithm}             &$\A~~|~~\B~~|~~\C$\\
\tabucline[1.2pt]\\
EXTRA~\cite{shi2015extra} & $\frac{1}{2}(\I+\W)~~|~~\I~~|~~\frac{1}{2}(\I-\W)$\\
NIDS~\cite{li2017decentralized}/Exact Diffusion~\cite{yuan2018exact_p1}          & $\frac{1}{2}(\I+\W)~~|~~\frac{1}{2}(\I+\W)~~|~~\frac{1}{2}(\I-\W)$ \\
NEXT~\cite{di2016next}/AugDGM~\cite{xu2015augmented}   & $\W^2~~|~~\W^2~~|~~(\I-\W)^2$  {} \\
DIGing~\cite{nedich2016achieving}/ \cite{qu2017harnessing}     & $\W^2~~|~~\I~~|~~(\I-\W)^2$  {} \\
\cite{jakovetic2018unification}  & $b\W^2+(1-b)\W~~|~~\I~~|~~b\W^2-(1+b)\W+\I$  {}\\
 \cite{mansoori2019general}     &$\W^K~~|~~\sum_{i=1}^{K-1}{\W^i}~~|~~\W-\W^K$  {}\\
\hline
 \cite{alghunaim2019linearly} ($G\neq 0$)      &$\W~~|~~\I~~|~~\alpha(\I-\W)$\vspace{0.1cm}\\
\end{tabu}\caption{Connections with existing distributed algorithms. All the schemes but ours and \cite{alghunaim2019linearly}  apply only to \eqref{prob:dop_nonsmooth_same} with $G=0$.}
\label{tab:connections_to_existing_algs}
\vspace{-0.8cm}
\end{table}

\section{Convergence Analysis}\vspace{-0.1cm}
We establish  linear rate  of Algorithm \eqref{alg:g-pd-ATC}
  under the following assumption (along with Assumption \ref{assum:cond_A_B}).

\begin{assum}\label{assum:conditions_convergence}
The weight matrices $A\in \mathbb{R}^{m\times m},\,  {B\in \mathbb{S}^{m}}$ and $C\in \mathbb{S}_+^{m} $ satisfy: i) 
  $A = BD$ for some $-I \prec D \preceq I$; ii) $0 \prec I-C$; iii)
$B$ and $C$ commute; and iv)
$B^2 \prec  \frac{(L+\mu)^2}{\left( L\lambda_{\text{max}}(D)-  \mu\lambda_{\text{min}}(D) \right)^2} (I-C) $.
\vspace{-0.2cm}
\end{assum}
 
 Assumption~\ref{assum:conditions_convergence} together with Assumption 3 are quite mild and satisfied by a variety of algorithms; for instance, this is the case for all the schemes  in Table I (see \cite{XuSunScutariJ} for more details). In particular, the commuting property is trivially satisfied when $B, C\in P_K(W)$, for some given $W\in \mathcal{W}_{\mathcal G}$ (as in Table I). Also, one can  show that condition iv) is {\it necessary} to achieve linear rate.
\begin{thm}
\label{thm:contraction_T_c_T_f}
 Consider Problem~\eqref{prob:dop_nonsmooth_same} under Assumption~\ref{assum:Lipschiz_gradient}, whose optimal solution  is $x^\star$. Let $\{(\x^k,\z^k,\y^k)\}_{k\geq 0}$ be the sequence generated by Algorithm \eqref{alg:g-pd-ATC} under Assumptions   \ref{assum:cond_C} and \ref{assum:cond_A_B} and step-size   \vspace{-0.3cm}  
\begin{align*}
& \frac{1}{\mu} \left( \lambda_{\max}(D) - \lambda_{\max}\left( B^2(I-C)^{-1}\right)^{- {1}/{2}}\right)_+ <\gamma  \\
& \quad < \frac{1}{L} \left( \lambda_{\min}(D) + \lambda_{\max}\left( B^2(I-C)^{-1}\right)^{- {1}/{2}}\right).
\end{align*}
Then $\norm{\x^k-1x^{\star\top}}^2=\mathcal{O}(\lambda^k)$, with   
\begin{equation} \label{def_lambda}
  \lambda \triangleq  \max \left(q^2 \lambda_{\max}(B^2 (I-C)^{-1}), ~1-\lambda_2(C) \right)<1, \vspace{-0.3cm}
\end{equation}
and \vspace{-0.2cm}
\begin{align}\label{eq:def_q}
q \triangleq & \max\left(\abs{\lambda_{\min}(D)-\gamma L}, ~ \abs{\lambda_{\max}(D)-\gamma \mu} \right) .
\end{align}
The optimal step-size is $\gamma^\star \triangleq \frac{\lambda_{\max}(D)+\lambda_{\min}(D)}{L+\mu}$ leading to  the  smallest $q^\star \triangleq \frac{L\lambda_{\max}(D)-  \mu\lambda_{\min}(D)}{L+\mu}$, and thus the optimal rate.
\end{thm}

 \begin{col}
\label{col:contraction_overall_optimal_choice}
Under the same setting as Theorem~\ref{thm:contraction_T_c_T_f},  let $B^2\preceq I-C$ and $\A=\B$, so that $D=I,~\gamma^\star =\frac{2}{L+\mu}$. Then, the rate reduces to \vspace{-0.1cm}
\begin{equation}\label{rate-separation}
\lambda=\max\left\{ \bracket{\frac{\kappa-1}{\kappa+1}}^2, ~ 1-\lambda_2(\C) \right\}.\vspace{-0.2cm}
\end{equation}
\end{col}
Note that the lower bound condition on the step-size in Theorem \ref{thm:contraction_T_c_T_f} nulls when  $B^2(I-C)^{-1}\preceq I$ (since $\lambda_{\max}(D)=1$). Theorem  \ref{thm:contraction_T_c_T_f} and Corollary \ref{col:contraction_overall_optimal_choice} provide a unified set of convergence conditions for CTA- and ATC-based distributed algorithms. We refer to \cite{XuSunScutariJ} for a detailed discussion of several special instances.  Here, we mainly comment   Algorithm \eqref{alg:g-pd-ATC}  in the setting of Corollary \ref{col:contraction_overall_optimal_choice}.
This special instance enjoys two desirable properties, namely: \textbf{(i) rate-separation:}  
The   rate (\ref{rate-separation}) is determined by the worst rate  between the one  due to the communication $[1-\lambda_2(\C)]$ and that of the optimization $[((\kappa-1)/(\kappa+1))^2]$. This separable structure is the key enabler   for our distributed scheme to achieve the convergence rate of the centralized proximal gradient algorithm applied to Problem~\eqref{prob:dop_nonsmooth_same}--see Sec.~\ref{sec:tradeoff}; and  \textbf{(ii) network-independent step-size:} The step-size in Corollary \ref{col:contraction_overall_optimal_choice} does not depend on the network parameters but only on the optimization and its value coincides with the optimal step-size of the centralized proximal-gradient algorithm. This is a major advantage over current distributed schemes applicable to \eqref{prob:dop_nonsmooth_same} (with $G\neq 0$) and complements the results in  \cite{li2017decentralized}, whose algorithm  however  cannot deal with  the non-smooth term  $G$ and use a non-optimal step-size.\vspace{-0.2cm}

\section{  Communication and computation trade-off}
\label{sec:tradeoff}
In this section we build on the rate separation property in Corollary 4   to show how to choose the matrices   $A,\,B$ and $C$ to  achieve the same rate of the centralized proximal gradient algorithm,  possibly using multiple (finite) rounds of communications.
\\\indent
Note  that    $\rho_{\texttt{opt}}\triangleq (\kappa-1)/(\kappa+1)$  is the  rate of the centralized proximal-gradient algorithm applied to  Problem~\eqref{prob:dop_nonsmooth_same}, under Assumption 1.
This means that if the network is ``well connected'', specifically   $1-\lambda_2(\C)\leq \rho_{\texttt{opt}}^2$, the proposed algorithm with the choice of  $A,\,B$ and $C$ under consideration already converges at  the \emph{desired}  linear rate $\rho_{\texttt{opt}}$. 
 On the other hand, when  
$1-\lambda_2(\C)>\rho_{\texttt{opt}}^2$, one can still achieve the centralized rate $\rho_{\texttt{opt}}$   by enabling multiple (finite) rounds of communications per proximal gradient evaluations. 
 We discuss next two strategies to reach this goal, namely: 1) performing multiple rounds of plain consensus using each time the same weight matrix; and 2) employing acceleration via Chebyshev polynomials. 
\\\noindent\textbf{1) Multiple rounds of consensus:}
Given a weight matrix $W\in \mathcal{W}_{\mathcal{G}}$ (i.e., compatible with $\mathcal{G}$),   we consider two possible choices of $\A,\B,\C$   satisfying Corollary \ref{col:contraction_overall_optimal_choice} and leading to distributed algorithms.  \textbf{Case 1:} Suppose $W\in\mathbb{S}^m_{++}$. We set  $\A=\B=\I-\C=W$, which implies $\B^2\preceq \I-\C$ (cf. Corollary~\ref{col:contraction_overall_optimal_choice}).  The resulting algorithm implemented using   (\ref{eq:g-pd-ATC_eliminate_y}) or (\ref{eq:g-pd-ATC_eliminate_y_G=0})  will require 
one communication exchange per gradient evaluation.  Note that this setting subsumes most existing primal-dual methods such as NIDS~\cite{li2017decentralized}/Exact Diffusion~\cite{yuan2018exact_p1}. If $W$ in the setting above is replaced by $\W^K$, with $K>1$, this corresponds to run 
$K$ rounds of consensus per computation, each round  using  $\W$.   Denote  
$\rho_{\texttt{com}}\triangleq\lambda_{\max}(\W-\J)$; we have $1-\lambda_2(\C)=\lambda_{\max}(W^K-\J)=\rho^K_{\texttt{com}}$. The value of $K$ is chosen to minimize the resulting rate $\lambda$ [cf. (\ref{rate-separation})], i.e., such that $\rho_{\texttt{com}}^K\leq\rho_{\texttt{opt}}^2$, which leads to $K=\lceil\log_{\rho_{\texttt{com}}}({\rho^2_{\texttt{opt}}})\rceil$. \textbf{Case 2:} Consider now the case  $\W\in\mathbb{S}^m$ and $\text{det}(W)\neq 0$. We can set  $\A^2=\B^2=\I-\C=W^2$, so that   Corollary~\ref{col:contraction_overall_optimal_choice} still applies. With this choice, every update in (\ref{eq:g-pd-ATC_eliminate_y}) or (\ref{eq:g-pd-ATC_eliminate_y_G=0}) will call for two communication exchanges per gradient evaluation. 
To reach the centralized   rate   $\rho_{\texttt{opt}}^2$, the optimal $K$ can be still found   as  $1-\lambda_2(\C)=(\lambda_{\max}(\A^{2K}-\J))=(\lambda_{\max}(\A-\J))^{{2K}}\leq  \rho_{\texttt{opt}}^2$.\\ 
\noindent \textbf{2) Chebyshev acceleration:}
To further reduce the number of communication steps, we can leverage Chebyshev acceleration~\cite{auzinger2011iterative}. Specifically, in the setting of Case 2 above, we set   $\A=\P_K(\W)$ and $P_K(1)=1$ (the latter is to ensure the double stochasticity of $\A$), with $P_K\in\mathbb{P}_K$. 
This leads to $1-\lambda_2(C)=\lambda_{\max}(\A^2-\J)$.  {The idea of Chebyshev acceleration is to find the ``optimal'' polynomial $P_K$ such that $\lambda_{\max}(\A^2-\J)$ is minimized,
 i.e.,
$\rho_C\triangleq\min_{P_K\in \mathbb{P}_K,P_K(1)=1}\max_{t\in[-\rho_{\texttt{com}},\rho_{\texttt{com}}]} \abs{P_K(t)}$.
The optimal solution of this problem  is  $P_K(x)={T_K(\frac{x}{\rho_{\texttt{com}}})}/{T_K(\frac{1}{\rho_{\texttt{com}}})}$    \cite[Theorem 6.2]{auzinger2011iterative}, with $\alpha'=-\rho_{\texttt{com}}$, $\beta'=\rho_{\texttt{com}},\gamma'=1$ (which are certain parameters therein),
where $T_K$ is the $K$-order Chebyshev polynomials that can be computed in a distributed manner via the following  iterates \cite{auzinger2011iterative,scaman17optimal}:
$T_{k+1}=2\xi T_k(\xi)-T_{k-1}(\xi),$  $k\geq 1$,
with $T_0(\xi)=1$, $T_1(\xi)=\xi$. 
  Also, invoking \cite[Corollary 6.3]{auzinger2011iterative}, we have
$\rho_C=\frac{2c^K}{1+c^{2K}}$
where $c=\frac{\sqrt{\vartheta}-1}{\sqrt{\vartheta}+1},\vartheta=\frac{1+\rho_{\texttt{com}}}{1-\rho_{\texttt{com}}}$. Thus, the minimum value of $K$ that leads to $\rho_C\leq\rho^2_{\texttt{opt}}$ can be obtained as $K=\lceil\log_c^{ 1/\rho^2_{\texttt{opt}}+\sqrt{1/\rho^4_{\texttt{opt}}-1}   }\rceil$. Note that to be used in the setting above, $A$ must be returned as nonsingular.  
\begin{figure*}[t]
\centering
\subfloat[\cite{van2019distributed}, $\A=\W^K$ ]{\includegraphics[width=0.3\textwidth]{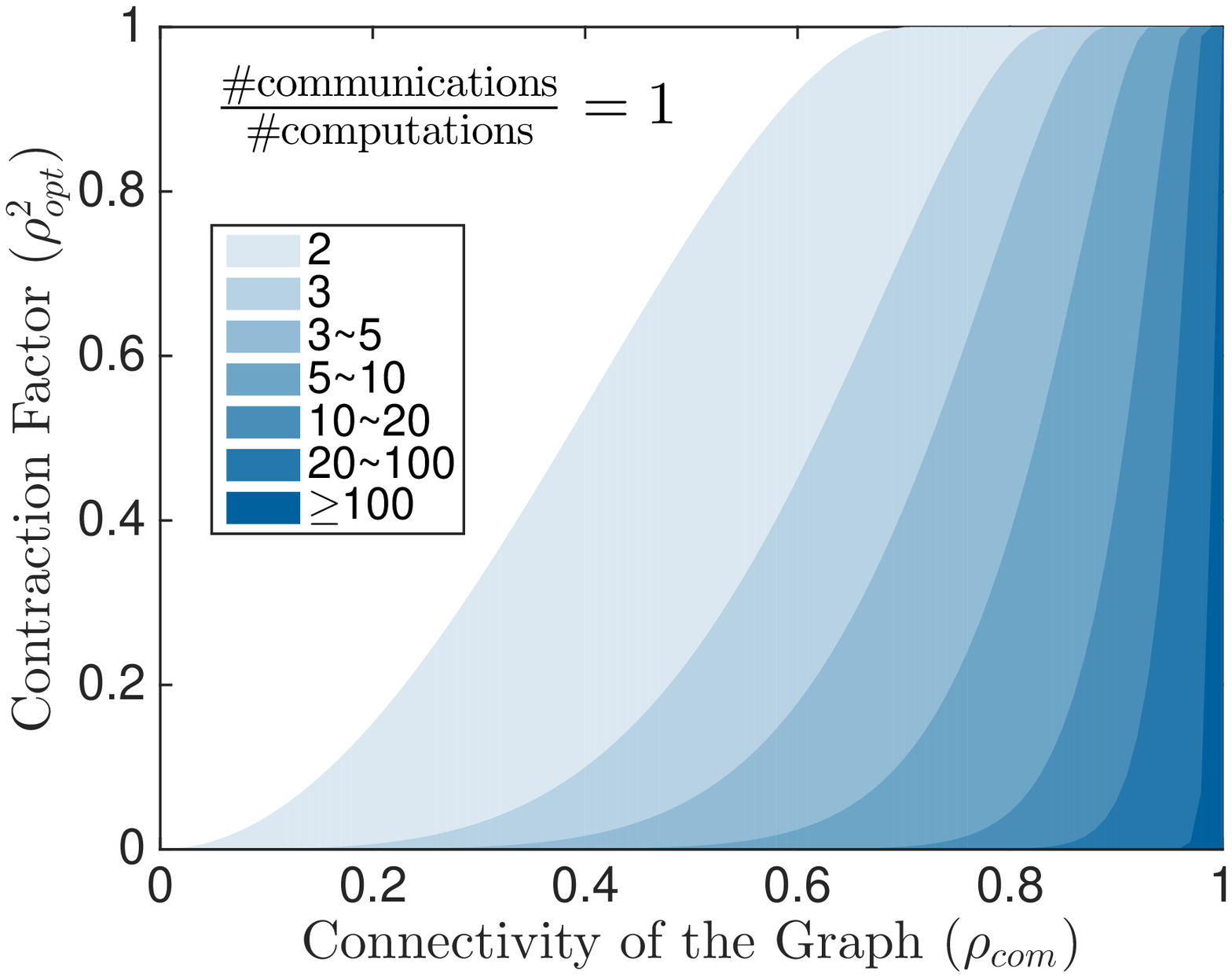}\label{fig_first_case}}
\hfil
\subfloat[Proposed scheme,   $\A=\W^K$]{\includegraphics[width=0.3\textwidth]{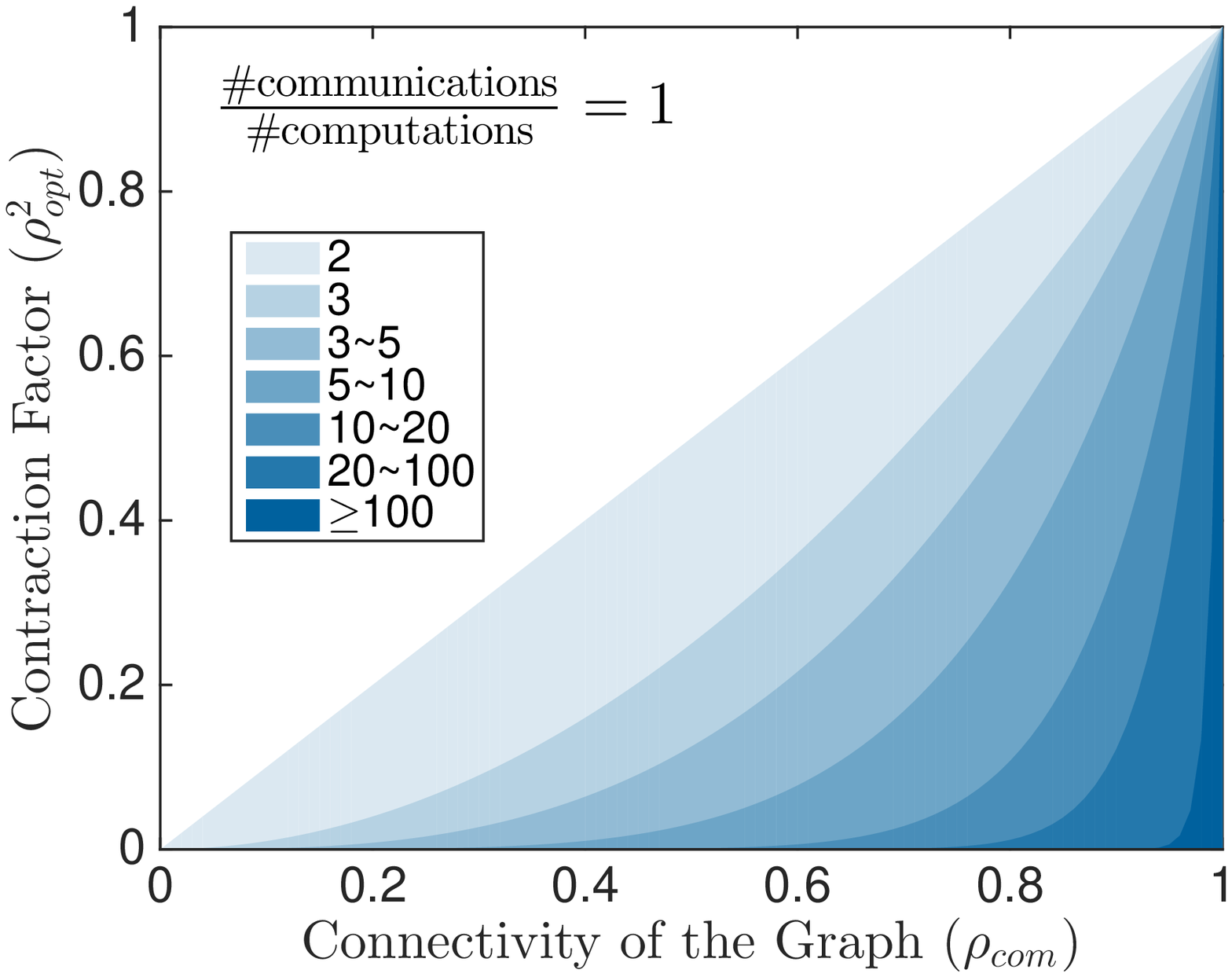}%
\label{fig_second_case}}
\hfil
\subfloat[Proposed scheme, $\A=P_K(\W)$]{\includegraphics[width=0.3\textwidth]{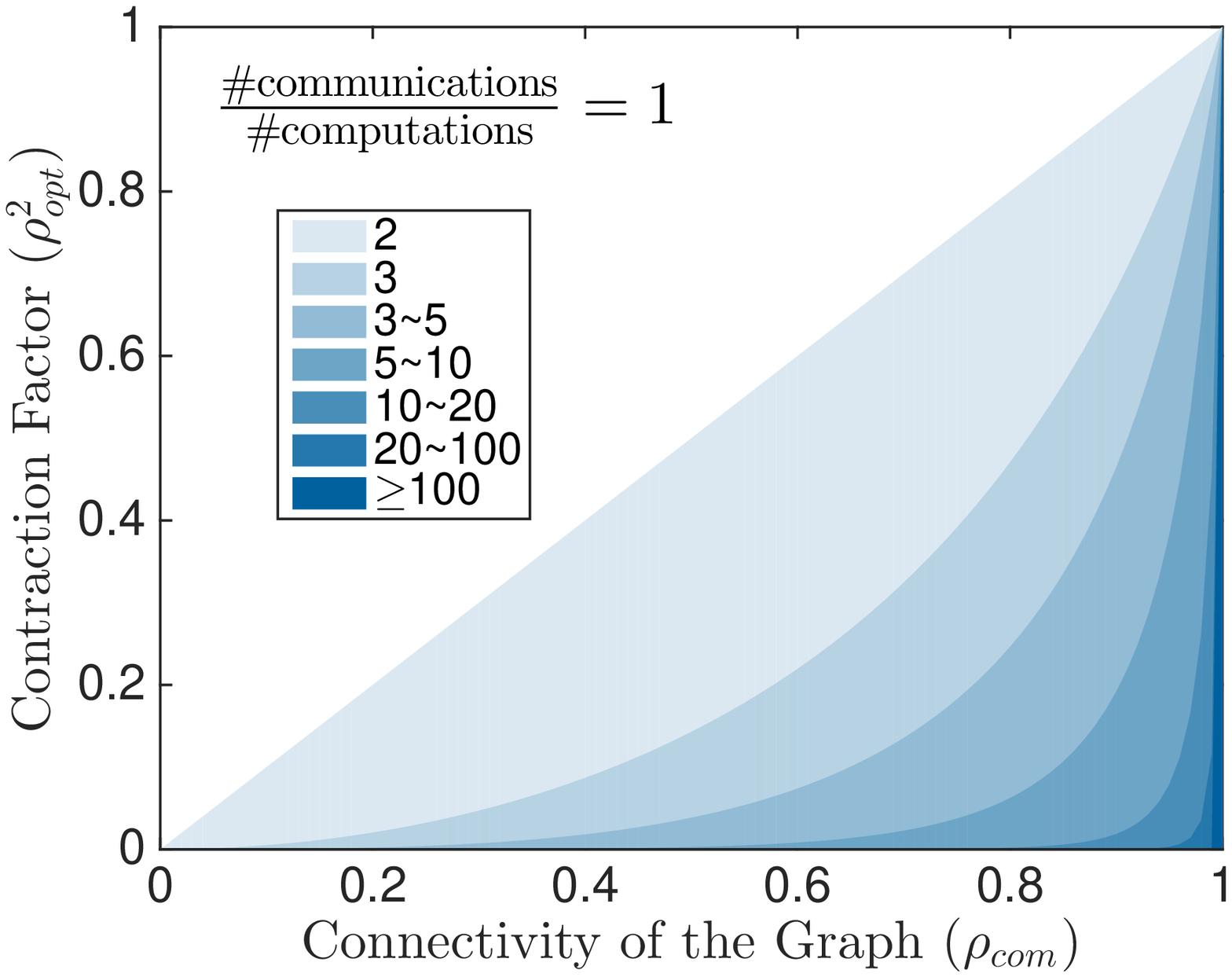}%
\label{fig_second_case}}\vspace{-0.1cm}
\caption{Ratio between the number of communications and
computations to achieve the centralized linear rate, as a function of the spectral gap $\rho_{\texttt{com}}$  and the gradient contraction factor $\rho_{\texttt{opt}}$. The proposed scheme employing multiple consensus rounds (subplot (b)) and Chebyshev acceleration (subplot (c)) is compared with \cite{van2019distributed} (subplot (a)).}\vspace{-0.4cm}
\label{fig_tradeoff}
\end{figure*}
\\\indent In Fig.~\ref{fig_tradeoff} we plot the minimum number $K$ of communication steps needed to achieve the rate of the centralized gradient   as a function of $\rho_{\texttt{com}}$ and $\rho^2_{\texttt{opt}}$. Since only one computation is performed per iteration, this adjusts the ratio between the number of communications and computations. We compare our algorithm in the setting of Case 2 above,  using $\A=\W^K$ or Chebyshev acceleration  $\A=P_K(\W)$, with the distributed scheme in  \cite{van2019distributed}.  The figure shows that    (i)     Chebyshev acceleration helps to reduce the number of communications   to sustain a given rate; and (ii) when $\rho_{\texttt{opt}}$ is close to $1$ ($\kappa$ is ``large''),  {both instances of the proposed scheme need much less communication steps to attain the centralized rate than that in \cite{van2019distributed}.  
 More specifically, to match the rate $\rho_{opt}$, one needs to run at least $K$ number of communications such that:\vspace{-0.1cm}
\[
\begin{aligned}
\rho_{com}^K=
\begin{cases}
\rho_{opt}^2,&\text{[this work]};\\
\frac{\sqrt{1+\rho_{opt}}-\sqrt{1-\rho_{opt}}}{2}, &\text{\cite{van2019distributed}}.
\end{cases}
\end{aligned}
\]
When $\rho_{opt}\rightarrow 0$, we have $({\sqrt{1+\rho_{opt}}-\sqrt{1-\rho_{opt}}})/{2}\approx {\rho_{opt}}/{2}$. Thus, $\rho_{opt}^2\leq {\rho_{opt}}/{2}$, since $\rho_{opt}\ll{1}/{2}$; hence, the scheme in \cite{van2019distributed}   needs less number of communications than the proposed algorithm in the aforementioned setting. On the other hand, when $\rho_{opt}\rightarrow 1$, we have $({\sqrt{1+\rho_{opt}}-\sqrt{1-\rho_{opt}}})/{2}\approx {\rho_{opt}}/{\sqrt{2}}$. In this case,   $  {\rho_{opt}}/{\sqrt{2}}\leq {1}/{\sqrt{2}}\leq \rho_{opt}^2$; hence, our scheme require  less  communications than that in \cite{van2019distributed}. Moreover, since $({\sqrt{1+\rho_{opt}}-\sqrt{1-\rho_{opt}}})/{2}\leq  {1}/{\sqrt{2}}<1$,  when $\rho_{com}\rightarrow 1$, the scheme in \cite{van2019distributed}  will  need significantly more communication    to match the centralized optimal rate.\vspace{-0.2cm}

\section{Conclusion}\vspace{-0.1cm}
 We proposed a unified distributed algorithmic framework for composite optimization problems over   networks; the algorithm includes many existing schemes as special cases. Linear rate was proved,  leveraging a  contraction operator-based anaysis.
Under a proper choice of the design parameters, the rate dependency on the network and cost functions can be decoupled, which allowed us to  determine the minimum number of communication steps needed to match the   rate of centralized (proximal)-gradient methods.\vspace{-0.2cm}

\appendix

\vspace{-0.1cm}

\section{Convergence Analysis}
We provide here a sketch of the proof of Theorem  \ref{thm:contraction_T_c_T_f}; see \cite{XuSunScutariJ} for more details.  Assumptions 2 and 3 are tacitly assumed hereafter.

\noindent {\bf Step 1: Auxiliary sequence and operator splitting:} Lemma \ref{lem-tranformation} below interpretes  \eqref{alg:g-pd-ATC}  as the fixed-point iterate of a suitably defined composition of   contractive and nonexpansive operators. 
\begin{lem}[\!\cite{XuSunScutariJ}]\label{lem-tranformation}
Given the sequence $\{(\z^k, \x^k, \y^k)\}_k$ generated by Algorithm \eqref{alg:g-pd-ATC}, define $U^k \triangleq 
[({\z}^k)^\top, ({\y}^k)^\top]^\top$. There holds
\[ U^{k}=
\begin{bmatrix}B & 0\\
0 & B\sqrt{C}
\end{bmatrix} 
\underset{\widetilde{U}^k}{\underbrace{\begin{bmatrix}\widetilde{\z}^{k}\\
\sqrt{C}\widetilde{\y}^{k}
\end{bmatrix}}},\vspace{-0.2cm}\]
with $\{\widetilde{U}^k\}_k$ defined by the following dynamics
\begin{align*}
\widetilde{U}^{k+1} = \underbrace{
\begin{bmatrix}
(D- \gamma \nabla f)\circ \text{prox}_{\gamma g} \circ B  & -\sqrt{C}\\
\sqrt{C}(D- \gamma \nabla f)\circ \text{prox}_{\gamma g} \circ B   & I-C
\end{bmatrix}}_{T}
\widetilde{U}^k,\quad k\geq 1,
\end{align*}
 and the initialization  $\widetilde{\z}^{1} = \widetilde{\y}^{1} = (D - \gamma \nabla f) (\x^{0})$. Furthermore, the operator $T$ can be decomposed as
\begin{align*}
T=
\underbrace{
\begin{bmatrix}
I  & -\sqrt{C}\\
\sqrt{C}   & I-C
\end{bmatrix}}_{\triangleq T_C}
\underbrace{
\begin{bmatrix}
D- \gamma \nabla f & 0\\
0  & I
\end{bmatrix}}_{\triangleq T_f}
\underbrace{
\begin{bmatrix}
\text{prox}_{\gamma g}  & 0 \\
0  & I
\end{bmatrix}}_{\triangleq T_g}
\underbrace{
\begin{bmatrix}
B &  0 \\
0    & I
\end{bmatrix}}_{\triangleq T_B},
\end{align*} where   $T_C$ and $T_B$ are the operators associated with communications while    $T_f$   and $T_g$ are the   gradient and proximal operators, respectively. Finally, every fixed point $\widetilde{U}^\star\triangleq [\widetilde{Z}^\star,\sqrt{C}\widetilde{Y}^\star]$ of $T$ is such  that $B \widetilde{Z}^\star=  1 x^{\star \top}\in \mathcal{S}_{\texttt{Fix}}$. \end{lem}  \vspace{-0.3cm}
Building on Lemma \ref{lem-tranformation}, the proof of Theorem  \ref{thm:contraction_T_c_T_f} reduces to showing  $\| \widetilde{\z}^k- \widetilde{\z}^\star\|=\mathcal{O}(\lambda^k)$. To do so, Step 2 below studies the contraction (nonexpansive) properties of single   operators composing $T$ while   Step 3 chains these properties showing that $T$ is  {$\lambda$-contractive with respect to a suitable norm}. 

\noindent \textbf{Step 2: On the properties of   $T_C$, $T_f$, $T_g$  and $T_B$.} We summarize next the main properties of the aforementioned operators; proofs of the results below can be found in \cite{XuSunScutariJ}. We will use the following notation: given  $X \in \mathbb{R}^{2m\times d}$, we denote by  $(X)_u$ and $(X)_\ell$  its upper and lower $m\times d$ matrix-block.   
\begin{lem} The operator $\T_c$ satisfies
\label{lem:contraction_T_c}
\[
\norm{\T_C\, X-\T_C \,Y }^2_{\Lambda_C} = \norm{X-Y }^2_{V_C},\quad \forall X,Y\in\mathbb{R}^{2m\times d},
\]
where $\Lambda_C \triangleq  \diag(\I-C, \I)$ and $V_C\triangleq\diag(\I,\I-C)$.
\end{lem}

\begin{lem}\label{lem:contraction_T_f}
 With $q$ defined in in Th.  \ref{thm:contraction_T_c_T_f}, $\T_f$ satisfies:  $\forall X, Y \in\mathbb{R}^{2m\times d}$,
\begin{equation*}
  \norm{(\T_f \,X)_u - (\T_f\,Y)_u}^2\leq q^2 \norm{(X)_u- (Y)_u}^2\,\,\text{and}\,\,  (\T_f\,X)_\ell = (X)_\ell.
\end{equation*}    
\end{lem}
\begin{lem} $\T_g$ satisfies:  $\forall X, Y\in\mathbb{R}^{2m\times d}$,
\label{lem:contraction_T_g}
\begin{equation*}
  \norm{(\T_g\, X)_u-(\T_g\, Y)_u}^2\leq \norm{(X)_u-(Y)_u}^2 \,\,\text{and}\,\,
  (\T_g\,X)_\ell = (X)_\ell.
\end{equation*}
\end{lem}

\begin{lem} The operator  $\T_B$ satisfies:
\label{lem:contraction_T_B}
\begin{equation*}
\norm{(\T_B\, X)_u}^2 = \norm{(X)_u}_{B^2}^2,\quad (\T_g\,X)_\ell = (X)_\ell,\quad \forall X \in\mathbb{R}^{2m\times d}.
\end{equation*}
\end{lem}
\vspace{-0.2cm}
\noindent \textbf{Step 3: Chaining Lemmata 6-9.}
 Define the matrices $Q_f \triangleq \diag(q^2 \I,\,\I)$ and $\Lambda_B=\diag(B^2, \I);$
the  contraction property of  $T$ are implied by the following chain: $\forall \, X, Y\in\mathbb{R}^{2m\times d}$ with $X_\ell,Y_\ell \in \Span{\sqrt{C}}$,
\begin{align*}
& \norm{T\,X- T\,Y}_{\Lambda_C}^2   
 \stackrel{Lm.~\ref{lem:contraction_T_c}}{=}\norm{T_f \circ T_g \circ T_B\,(X- Y)}_{V_C}^2 \\
& \stackrel{Lm.~\ref{lem:contraction_T_f}}{\leq} \norm{ T_g \circ T_B\,(X- Y)}_{V_C Q_f}^2  
  \stackrel{Lm.~\ref{lem:contraction_T_g}}{\leq} \norm{ T_B\,(X- Y)}_{V_C Q_f}^2 \\
& \stackrel{Lm.~\ref{lem:contraction_T_B}}{=} \norm{X- Y}_{V_C Q_f \Lambda_B}^2 \stackrel{(*)}{\leq} \lambda \, \norm{X- Y}_{\Lambda_C}^2,
\end{align*}
where $V_C Q_f \Lambda_B = \diag(q^2 B^2, I-C)$,   $\lambda$ is defined in   (\ref{def_lambda}); and    (*) is due to the following two facts:  i)   
 $ q^2 \|(Z)_u\|^2_{B^2}=q^2 \|(I-C)^{\frac{1}{2}}(Z)_u\|^2_{B^2(I-C)^{-1}}   \leq q^2 \lambda_{\text{max}} (B^2(I-C)^{-1})\|(I-C)^{\frac{1}{2}}(Z)_u\|^2 
  = q^2 \lambda_{\text{max}} (B^2(I-C)^{-1})\norm{(Z)_u}^2_{I-C},$ for all $\,(Z)_u \in \mathbb{R}^{m\times d}$; 
 and ii) $X_\ell,Y_\ell \in \Span{\sqrt{C}}$.

    \bibliographystyle{IEEEtran}
\bibliography{IEEEabrv,reference}
 
 \end{document}